\documentclass[12pt,a4paper]{article}
\pdfoutput=1
\usepackage[italian]{babel}
\usepackage{amsmath}
\usepackage{amssymb}
\usepackage{latexsym}
\usepackage{verbatim}
\usepackage{graphicx}
\DeclareGraphicsExtensions{.png}

\begin{document}
\begin{center}
\Large
\textbf{La cometa di Goldbach e ... le \vspace{6mm} altre} \\
\normalsize
\textbf{Donato\; SAELI,\; \ Maurizio\; \vspace{8mm} SPANO}
\begin{minipage}{120mm}
\small
\textsc{Abstract.} Goldbach's comet is the plot of the Goldbach function $\, g(n), \,$ in the interval $\, [3,N], \,$ with a large positive integer $\, N. \,$ The function $\, g(n) \,$ counts the number of different ways in which $\, 2n \,$ can be expressed as the sum of two odd primes. An account, hopefully satisfying and accessible, is given for the layers that make up the comet. By means of several (sometimes historical) results of Theory of Number, other conjectures, similar to the Goldbach's one, emerge. These are related with sequences of odd positive integers, like but not quite to the prime \vspace{2mm} sequence. \\
\textsc{Keywords:} Goldbach's comet, function, conjecture, extended conjecture, Sylvester factor. Prime number theorem, PNT for arithmetic progressions, asymptotic expression for the nth \vspace{2mm} prime. \\
\textsc{MSC:} \vspace{6mm} 11P32
\end{minipage}
\end{center}
\linespread{1.3}
\normalsize
\paragraph{1 La cometa di Goldbach.}
Il grafico della funzione aritmetica $\, g(n) \,$ che associa ad $\, n \,$ il numero delle coppie $\, (p,q) \,$ di numeri primi dispari tali che $\, p + q = 2n, \,$ considerata nell'intervallo $\, [3,N], \,$ con $\, N \,$ sufficientemente grande, appare curiosamente come una cometa (fig.$\, 1$) e prende giustappunto il nome di ``\textit{cometa di Goldbach}\,'';
\begin{center}
\includegraphics[scale =.4]{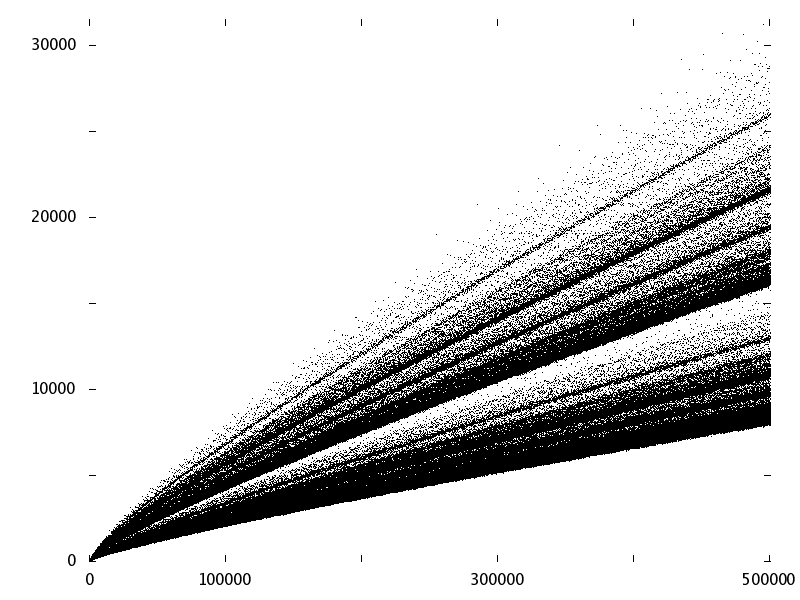} \\
\footnotesize Figura 1. Grafico di $\, g(n), \,$ per $\, n \in [\, 3 \, , \, 500.000\, ] \,$
\normalsize
\end{center}
tale nome deriva dal fatto che l'affermazione \ ``$g(n) > 0 \ \ per \ ogni \ \ n > 2 \,$'' \\
\`e equivalente alla congettura di Goldbach: \\
``\textit{Ogni numero pari non inferiore a quattro \`e somma di due primi}\,''. \footnote{ \ \ La congettura risale al 1742 (cfr. [D], p. 421) e appare tuttora aperta, sebbene \\
\hspace*{7.5mm} parecchi studiosi abbiano conseguito risultati considerevoli (cfr. ad es. [V], \vspace{1mm} [E], \\
\hspace*{7.5mm} [C] e [MV]).} \\
Un particolare che si nota in fig.$\, 1$ \`e la ``striatura della cometa\,''; si distinguono chiaramente due zone principali a loro volta divise in pi\`u striscie. \\
Di questa particolarit\`a del grafico della $\, g(n), \,$ vogliamo dare una giustificazione che speriamo sia soddisfacente e accessibile. \\
Naturalmente pietra miliare di ogni discussione su questi argomenti rimane la  ``Congettura A'', detta anche congettura ``\textit{estesa}\,'' di Goldbach [HL], formulata nel 1922 da Hardy e Littlewood:
$$g(n) \sim h(n) = \dfrac{4cn}{(\lg n)^2} \prod_{
\begin{subarray}{c}
p \, | \, n \\
p \geq 3
\end{subarray}}
\dfrac{p-1}{p-2} \; ,$$
\hspace*{15mm} dove \vspace{2mm} $\, \displaystyle c = \prod_{p \geq 3} \dfrac{p(p-2)}{(p-1)^2}$.\, \footnote{\ \ La scrittura $\, g(n) \sim h(n) \,$ indica che $\, \displaystyle \lim_{n\rightarrow\infty}\dfrac{g(n)}{h(n)} = 1 \,$ e si dice che le due \vspace{1mm} funzioni \\
\hspace*{7.5mm} sono \textit{asintoticamente equivalenti} \vspace{2mm} (per $\, n \rightarrow \infty$). \\
\hspace*{7.5mm} Il prodotto $\, \displaystyle \prod_{
\begin{subarray}{c}
p \, | \, n \\
p \geq 3
\end{subarray}}
\dfrac{p-1}{p-2} \,$ s'intende esteso a tutti i numeri primi dispari che dividono $\, n$ \\
\hspace*{7.5mm} e si pone uguale a 1 se \`e privo di fattori, cio\`e se \vspace{2mm} $\, n = 2^k.$ \\
\hspace*{7.5mm} L'espressione che determina la costante $\, c \,$ \`e un ``\textit{prodotto infinito}\,'' esteso a \underline{tutti} \\
\hspace*{7.5mm} i primi dispari, \ pi\`u precisamente:
$$c = \prod_{p \geq 3} \dfrac{p(p-2)}{(p-1)^2} = \lim_{n\rightarrow\infty} \bigg(\prod_{3 \leq p \leq n} \dfrac{p(p-2)}{(p-1)^2}\bigg) = 0,6601618 \dots$$
\hspace*{7.5mm} Un risultato prossimo alla congettura A era stato annunciato da Sylvester [S] \\
\hspace*{7.5mm} nel 1871, che aveva proposto una formula equivalente alla $\, g(n) \sim 2e^{-\gamma}h(n);$ \\
\hspace*{7.5mm} ma la congettura A ``...  \`e la sola formula di questa sorta che pu\`o essere corretta, \\
\hspace*{7.5mm} cosicch\'e la formula di Sylvester \`e errata. Ma Sylvester \`e stato il primo ad identificare \\
\hspace*{7.5mm} il \vspace{-2mm} fattore
$$\prod_{
\begin{subarray}{c}
p \, | \, n \\
p \geq 3
\end{subarray}}
\dfrac{p-1}{p-2}$$
\hspace*{7.5mm} a cui sono dovute le \textit{irregolarit\`a} della $\, h(n). \,$ Non vi sono indicazioni sufficienti \\
\hspace*{7.5mm} per mostrare come sia stato condotto al suo risultato. ...'' ([HL], pp. 32, \vspace{1mm} 33).}
\begin{center}
\includegraphics[scale =.4]{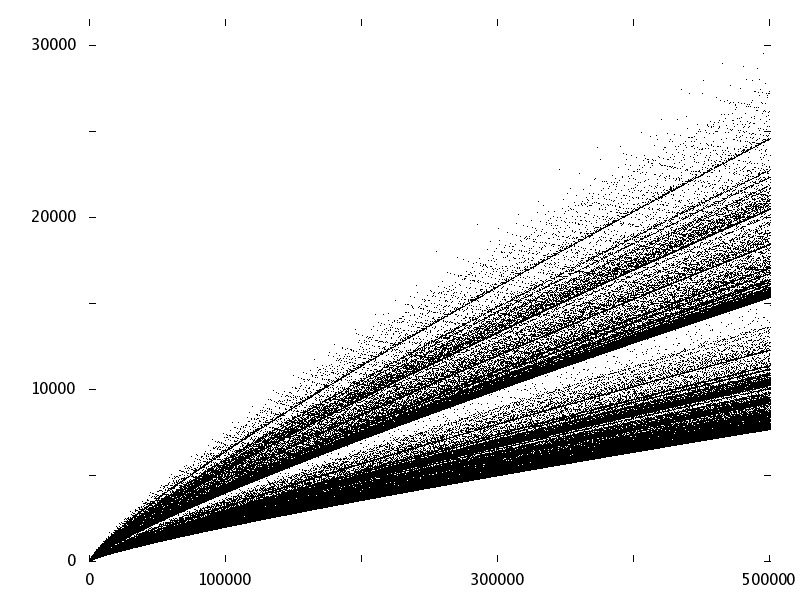} \\
\footnotesize Figura 2. Grafico di $\, h(n), \,$ per $\, n \in [\, 3 \, , \, 500.000\, ] \,$
\normalsize
\end{center}
La funzione $\, h(n), \,$ il cui grafico (fig.$\, 2$) somiglia alla cometa di Goldbach \\
(ne pare la ``bella copia\,''), \`e prodotto della funzione ``tranquilla'' \vspace{1.5mm} $\, \dfrac{4cn}{(\lg n)^2} \,$ (che, come si pu\`o vedere, limita il bordo inferiore delle comete) per il \vspace{1.5mm} fattore $\, \displaystyle \prod_{
\begin{subarray}{c}
p \, | \, n \\
p \geq 3
\end{subarray}}
\dfrac{p-1}{p-2} \,$ che dipende piuttosto ``vivacemente'' da $\, n \,$ e pu\`o essere relativamente grande se $\, n \,$ ha parecchi fattori primi piccoli; \footnote{ \ \ Il prodotto infinito $\, \displaystyle \prod_{p\geq 3} \dfrac{p-1}{p-2} \,$ esteso a tutti i primi dispari diverge a $\, + \infty$.} \ ad esempio vale $\, 3, \overline{5} \,$ per $\, n = 2310 = 2 \cdot 3 \cdot 5 \cdot 7 \cdot 11, \,$ ma meno di 1,0004331 per il numero primo $\, n = 2311. \,$ \quad D'altra parte \`e possibile mostrare l'importanza di questo fattore nella congettura A, contraendo la cometa di Goldbach in una stretta \vspace{1.5mm} scia; basta considerare la funzione $\, G(n) = g(n) \displaystyle \prod_{
\begin{subarray}{c}
p \, | \, n \\
p \geq 3
\end{subarray}}
\dfrac{p-2}{p-1} \,$  in luogo della $\, g(n), \,$ cosicch\'e la congettura A assume la forma:
$$G(n) \sim \dfrac{4cn}{(\lg n)^2}$$
ed infatti (fig.$\, 3$) il grafico della $\, G(n), \,$ la ``scia,'' si ``adagia'' sul \vspace{2mm} grafico \\
della $\,\dfrac{4cn}{(\lg n)^2}$.
\begin{center}
\includegraphics[scale =.4]{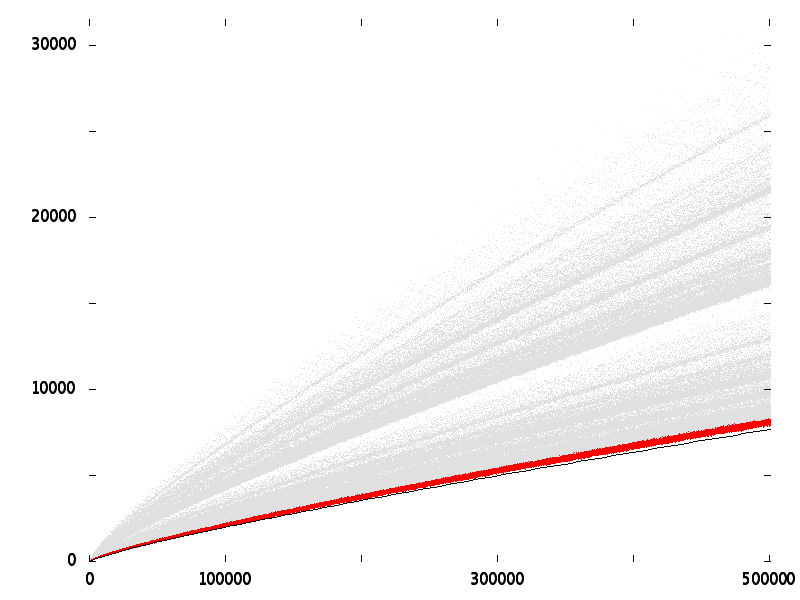} \\
\footnotesize Figura 3. Grafici delle funzioni $\, G(n) \ \ \textrm{e} \ \ \dfrac{4cn}{(\lg n)^2} \,$ per $\, n \in [\, 3 \, , \, 500.000\, ] \,$
\normalsize
\end{center}
\vspace{3mm}
\paragraph{2 La cometa di Goldbach a colori.}
Vi sono due aspetti, propri della distribuzione dei numeri primi fra i numeri naturali, che incidono profondamente sulla natura della funzione $\, g(n), \,$ \\
Il \textit{teorema dei numeri primi} (Hadamard, de la Vall\'ee-Poussin):
$$\quad \displaystyle \lim_{x\rightarrow +\infty}\dfrac{\pi(x)}{\dfrac{x}{\lg x}} = 1,$$
dove $\, \pi(x) \,$ associa a $\, x \,$ il numero dei primi che non superano $\, x; \,$ \\
l'altro consiste nell'estensione data da de la Vall\'ee-Poussin al teorema di Dirichlet sui \textit{primi in progressione aritmetica}: \\
\textit{Se\, m\, ed\, a\, sono due numeri naturali primi fra loro, allora}
$$\quad \displaystyle \lim_{x\rightarrow +\infty}\dfrac{\pi_{m,a}(x)}{ \quad \dfrac{x}{\varphi(m)\lg x} \quad } = 1,$$
\textit{dove} $\, \pi_{m,a}(x) \,$ \textit{associa a\, x\, il numero dei primi della forma\, a$+$km\, che non superano\, x \ e} \ $\, \varphi(m) = \raisebox{-1mm}{\includegraphics[scale =.11]{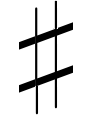}} \{h \in \mathbb{N}: 1 \leq h < m \ \textrm{e} \ (h,m) = 1 \} \,$ \textit{\`e ``l'indicatore\,'' di} Eulero; \ \textit{se} $\, (m,a) > 1, \,$ \textit{fra i numeri della forma\, a$+$km,\, uno al pi\'u pu\`o essere primo.} \\
Naturalmente \`e possibile e conveniente riformulare quest'ultimo teorema \\
nei termini dell'aritmetica modulare: \\
\textit{Nell'anello} $\, \mathbb{Z}_m \,$ \textit{delle classi di resto modulo\, m,\, una classe} $\, \bar{a} \,$ \textit{contiene infiniti numeri primi se e solo se \`e una classe prima col modulo\, m,\, cio\`e se e solo se} $\, (a,m) = 1; \,$ \textit{inoltre, se} $\, \bar{a} \ \ e \ \ \bar{b} \,$ \textit{sono \vspace{1mm} due classi prime col modulo\, m,\, allora} $\, \ \displaystyle \lim_{x\rightarrow +\infty}\dfrac{\pi_{m,a}(x)}{\pi_{m,b}(x)} = 1. \ \,$ \textit{Infine una \vspace{1mm} classe \underline{non} prima col modulo\, m\, pu\`o contenere al pi\`u un solo numero primo.} \\
Se $\, \bar{a} \in \mathbb{Z}_m, \,$ si dice \textit{``divisore''} della classe $\, \bar{a} \,$ e si indica con\, div($\bar{a}$)\, il massimo comun divisore $\, (a,m) \,$ fra $\, a \,$ e il modulo $\, m; \,$ si noti che se $\, a' \in \bar{a} \,$ \`e $\, (a',m) = (a,m). \,$ Cos\`{\i} $\, \bar{a} \,$ contiene infiniti numeri primi se e solo se il suo divisore \`e 1. \\
\`E chiaro che se il modulo $\, m \,$ \`e pari e $\, \bar{a} \,$ e $\, \bar{b} \,$ sono due classi prime con $\, m, \,$ allora  la classe $\, \bar{a} + \bar{b} = \overline{a+b} \,$ contiene solo numeri pari, vale a dire \`e una ``\textit{classe pari}\,''; recentemente [Dm] \`e stato mostrato che ogni classe p\underline{ari} si pu\`o esprimere come somma di due classi p\underline{rime} e anche come determinare il numero di tali rappresentazioni. \ Pi\`u precisamente ([S], pp. 2261-2263), indicato con $\, \mathbb{Z}_m^* \,$ l'insieme (\textit{gruppo}) delle classi prime con $\, m, \,$ se $\, \bar{c} \in \mathbb{Z}_m \,$ \\
(in generale) si ha:
$$\sigma_m(c) = \raisebox{-1mm}{\includegraphics[scale =.11]{dies}} \{(\bar{a},\bar{b}) \in (\mathbb{Z}_m^*)^2 :\bar{a}+\bar{b} = \bar{c} \, \} = m \prod_{
\begin{subarray}{c}
p \, | \, m \\
p \, | \, c
\end{subarray}
}\bigg(1-\dfrac{1}{p} \bigg) \prod_{
\begin{subarray}{c}
p \, | \, m \\
p \, \nmid \, c
\end{subarray}
} \bigg(1-\dfrac{2}{p} \bigg),$$
ma se si sceglie come modulo il prodotto dei primi $\, h \,$ numeri primi, \\
cio\`e $\, m = 2 \cdot 3 \cdots p_h \,$ \ e \ se \ $d = \textrm{div}(\bar{c}) = (c,m), \,$ si vede facilmente che
$$\sigma_m(c) =  \sigma_m(d) =  m \prod_{p \, | \, d}\bigg(1-\dfrac{1}{p} \bigg) \prod_{p \, | \, \frac{m}{d}} \bigg(1-\dfrac{2}{p} \bigg); \hspace{24mm} (1)$$
cos\`{\i}, per questo tipo particolare di modulo, $\, \sigma_m(c) \,$ dipende \underline{solo} dal divisore della classe cui appartiene $\, c. \,$ Ma l'insieme dei possibili divisori per le classi \underline{pari} di $\, \mathbb{Z}_m, \,$ cio\`e l'insieme dei divisori pari (e positivi) $\, d \,$ di $\, m, \,$ \`e costituito dai mumeri $\, \displaystyle d = 2 \cdot 3^{\alpha_2} \cdots p_h^{\alpha_h}, \,$ per tutte le possibili scelte di $\, \alpha_i \,$ su 0 o 1, per $\, i = 2, \dots , h; \,$ ne segue \vspace{-1.5mm} che
$$\sigma_m(d) = \sigma_m(2 \cdot 3^{\alpha_2} \cdots p_h^{\alpha_h}) = \prod_{i=2}^h (p_i-1 )^{\alpha_i} (p_i-2)^{1-\alpha_i}$$
in \vspace{1.5mm} particolare $\quad \displaystyle \sigma_m(2) = \prod_{i=2}^h (p_i-2) \, ; \ $ si noti che 2 \`e l'unico divisore di $\, m \,$ per cui $\, \sigma_m \,$ risulta \underline{dispari}, inoltre, per $\, d \,$ pari, \`e $\, \sigma_m(2) \leq \sigma_m(d).$ \\
Indicata con $\, \{p_i \}_{i \in \mathbb{N}} \,$ la successione dei numeri primi in ordine crescente $\, (p_1 = 2), \,$ sia $\, N \,$ un numero naturale sufficientemente grande; se $\, k \,$ \`e l'indice tale che $\, p_k \leq N < p_{k+1}, \,$ si assuma come modulo $\, m = 2 \cdot 3 \cdots p_k. \,$ \\
Si consideri ora un qualunque numero naturale $\, n, \ 3 \leq n \leq N, \,$ si ha:
$$d = (2n, \, m) = 2 \cdot 3^{\beta_2} \cdots p_k^{\beta_k}, \ \ \textrm{con} \ \ \beta_i = \left\{ \begin{array}{l}
1 \ \; \textrm{se} \ \; p_i | n \\
0 \ \; \textrm{se} \ \; p_i \nmid n  \\
\end{array} \right. \ \ \textrm{per} \ \ i = 2, \, \dots, \, k,$$
$\; \displaystyle \sigma_m(2n) = \sigma_m(d) = \prod_{i=2}^k (p_i-1 )^{\beta_i} (p_i-2)^{1-\beta_i} \, ; \qquad$ e finalmente si trova
$$\dfrac{\sigma_m(2n)}{\sigma_m(2)} = \dfrac{\displaystyle \prod_{i=2}^k (p_i-1 )^{\beta_i} (p_i-2)^{1-\beta_i}}{\displaystyle \prod_{i=2}^k (p_i-2)} = \prod_{i=2}^k \bigg( \dfrac{p_i-1}{p_i-2} \bigg)^{\beta_i} = \prod_{ \begin{subarray}{c}
p \, | \, n \\
p \geq 3
\end{subarray}}
\dfrac{p-1}{p-2} \, ,$$
il fattore di \vspace{2mm} Sylvester. \\
Per vedere cosa c'entra tutto ci\`o con i colori della cometa di Goldbach occorre esaminare alcuni casi relativi a valori particolari (successivi) del modulo $\, m.$ \\
\textbf{\textit{i}}) \ Per $\, m = 6 \,$ i numeri pari si distribuiscono nelle classi (pari) $\, \bar{0}, \ \bar{2} \ \, \textrm{e} \ \, \bar{4} \,$ con $\, \textrm{div}(\bar{0}) = 6 \ \, \textrm{e} \ \, \textrm{div}(\bar{2}) = \textrm{div}(\bar{4}) = 2, \,$ mentre i numeri primi, tranne 2 e 3, \\
si vanno a collocare nelle due classi prime $\, \bar{1} \ \, \textrm{e} \ \, \bar{5}. \,$ 
$$\textrm{Si ha:} \hspace{9mm} \bar{0} = \bar{1} + \bar{5} = \bar{5} + \bar{1}, \quad \bar{2} = \bar{1} + \bar{1} \quad \textrm{e} \quad \bar{4} = \bar{5} + \bar{5};$$
in accordo con la (1), ritenendo distinte le due espressioni $\, \bar{1} + \bar{5} \,$ e $\, \bar{5} + \bar{1} \,$ per la classe $\, \bar{0}; \,$ questa distinzione torna utile se si valutano le possibilit\`a di esprimere un determinato pari come somma di due primi appartenenti a una stessa classe oppure a due classi diverse. \hspace{\stretch{1}} Infatti, fissato $\, h \in \mathbb{N}, \,$ si \vspace{1mm} ha: \\
$6h-2 = 5+[6(h-2)+5] = 11+[6(h-3)+5] = 17+[6(h-4)+5] = \cdots = \\
\hspace*{17mm} [6(h-4)+5]+17 = [6(h-3)+5]+11 = [6(h-2)+5]+ \vspace{1mm} 5,$ \\
$6h \hspace{6.5mm} = 7+[6(h-2)+5] = 13+[6(h-3)+5] = 19+[6(h-4)+5] = \cdots = \\
\hspace*{17mm} [6(h-3)+1]+17 = [6(h-2)+1]+11 = [6(h-1)+1]+ \vspace{1mm} 5$ \\
$6h+2 = 7+[6(h-1)+1] = 13+[6(h-2)+1] = 19+[6(h-3)+1] = \cdots = \\
\hspace*{17mm} [6(h-3)+1]+19 = [6(h-2)+1]+13 = [6(h-1)+1]+ \vspace{2mm} 7;$ \\
si nota facilmente che mentre i numeri pari della forma $\, 6h \,$ si esprimono \\
in $\, h-1 \,$ modi diversi come somma di due addendi appartenenti rispettivamente alle classi $\, \bar{1} \ \, \textrm{e} \ \, \bar{5}, \,$ i pari della forma $\, 6h-2, \ \; \vspace{1mm} ( \ 6h+2 \ ) \,$ si possono esprimere, a meno dell'ordine, solo in $\, \bigg\lfloor \dfrac{h}{2} \bigg\rfloor \,$ modi diversi come somma di due \vspace{1mm} addendi entrambi appartenenti necessariamente alla classe $\, \bar{5} \ \; ( \ \bar{1} \ ).$ \\
Cos\`{\i} i numeri pari appartenenti alla classe $\, \bar{0}, \,$ rispetto ai numeri pari delle classi $\, \bar{2} \ \, \textrm{e} \ \, \bar{4}, \,$ hanno circa il doppio di possibilit\`a di rappresentazioni diverse come somma di due primi dispari. \footnote{\ \ Naturalmente sono state trascurate le possibilit\`a $\, 6h-2 = 3+[6(h-1)+1] \,$ \\
\hspace*{7.5mm} e $\, 6h+2 = 3+[6(h-1)+5].$} \\
Queste considerazioni suggeriscono come ``\textit{colorare la cometa di Goldbach.}'' Occorre semplicemente segnare, nell'illustrazione del grafico della $\, g(n), \,$ ad esempio, in nero i punti corrispondenti ai valori di $\, n \,$ per i quali $\, (2n,6) = 2 \,$ e in rosso gli altri; il risultato (fig.$\, 4$) pone in evidenza ancora di pi\`u le due zone principali.
\begin{center}
\includegraphics[scale =.4]{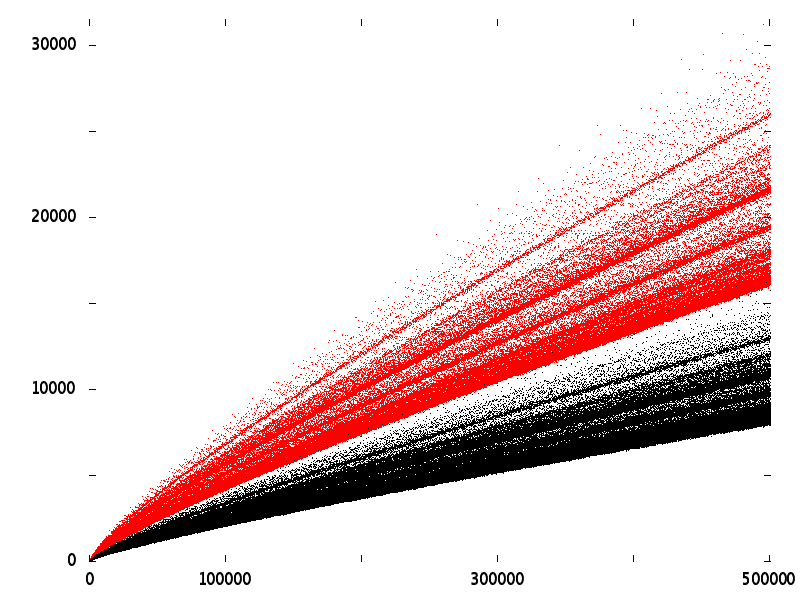} \\
\footnotesize Figura 4. Grafico di $\, g(n), \,$ per \vspace{-1.3mm} $\, n \in [\, 3 \, , \, 500.000\, ] \,$ \\
\hspace*{7mm} In nero \ i punti \vspace{-1.5mm} con $\, (2n,6) = 2,$ \\
\hspace*{7mm} in rosso i punti con $\, (2n,6) = 6.$ \\
\normalsize
\end{center}
\textbf{\textit{ii}}) \ Per $\, m = 30, \,$ vi sono otto classi prime: $\, \overline{1}, \, \overline{7}, \, \overline{11}, \, \overline{13}, \, \overline{17}, \, \overline{19}, \, \overline{23}, \, \overline{29} \,$ e quattro possibili divisori $\, d \,$ per le classi pari: $\, 2, \, 6, \, 10, \, 30 \,$ cui corrispondono, per la $\, (1), \,$ i rispettivi valori $\, 3, \, 6, \, 4, \, 8 \,$ di $\, \sigma_{_{30}}(d). \,$ Si pu\`o vedere facilmente che le classi $\, \overline{2}, \, \overline{4}, \, \overline{8}, \, \overline{14}, \, \overline{16}, \, \overline{22}, \, \overline{26}, \, \overline{28} \,$ hanno tutte divisore $\, 2 \,$ \\
e per queste classi $\, \sigma_{_{30}}(2) = 3; \quad $ infatti: $ \ \overline{2} = \overline{1} + \overline{1} = \overline{13} + \overline{19} = \overline{19} + \overline{13}, \,$ \\
$\, \overline{4} = \overline{11} + \overline{23} = \overline{17} + \overline{17} = \overline{23} + \overline{11}, \ \dots, \ \overline{28} = \overline{11} + \overline{17} = \overline{17} + \overline{11} = \overline{29} + \overline{29}.$ \\
Le classi $\, \overline{6}, \, \overline{12}, \, \overline{18}, \, \overline{24} \,$ hanno divisore $\, 6 \,$ e si ha: $\, \sigma_{_{30}}(6) = 6. \,$ \\
Le classi $\, \overline{10}, \, \overline{20} \,$ hanno divisore $\, 10 \,$ e $\, \sigma_{_{30}}(10) = 4. \,$ \\
Infine la classe $\, \overline{0} \,$ ha divisore $\, 30 \,$ e $\, \sigma_{_{30}}(30) = 8. \,$ \\
Cos\`{\i} ai quattro divisori pari di $\, 30 \,$ corrispondono altrettante zone nel grafico della $\, g(n) \,$ (fig.$\, 5$), le due zone principali vengono divise ciascuna in due \vspace{10mm} strati.
\begin{center}
\includegraphics[scale =.4]{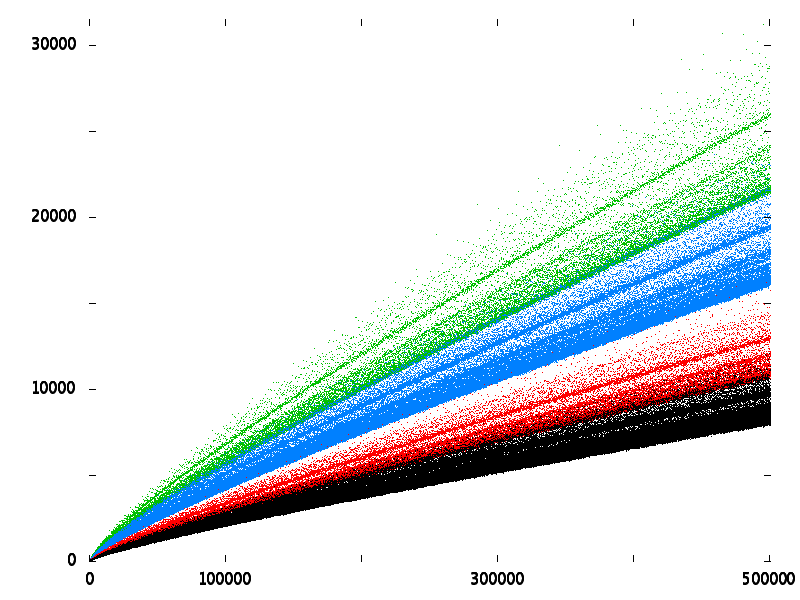} \\
\footnotesize Figura 5. Grafico di $\, g(n), \,$ per \vspace{-1.3mm} $\, n \in [\, 3 \, , \, 500.000\, ] \,$ \\
\hspace*{6mm} In nero\, i punti \vspace{-1.7mm} con $\, (2n,30) = 2,$ \\
\hspace*{7mm} in rosso i punti \vspace{-1.5mm} con $\, (2n,30) = 10,$ \\
\hspace*{8mm} in blu i punti \vspace{-1.7mm} con $\, (2n,30) = 6,$ \\
\hspace*{7mm} in verde i punti con \vspace{10mm} $\, (2n,30) = 30.$ \\
\normalsize
\end{center}
\textbf{\textit{iii}}) \ Per $\, m = 210, \,$ vi sono quarantotto classi prime e otto possibili \\
divisori $\, d \,$ per le classi pari, si ha:
\begin{center}
\begin{tabular}{c|cccccccc}
$d$ & 2 & 6 & 10 & 14 & 30 & 42 & 70 & 210 \\
\hline
$\sigma_m(d)$ & 15 & 30 & 20 & 18 & 40 & 36 & 24 & 48
\end{tabular}
\end{center}
ai divisori pari di $\, 30 \,$ corrispondono otto zone nel grafico della $\, g(n) \,$ (fig.$\, 6$).
\begin{center}
\includegraphics[scale =.4]{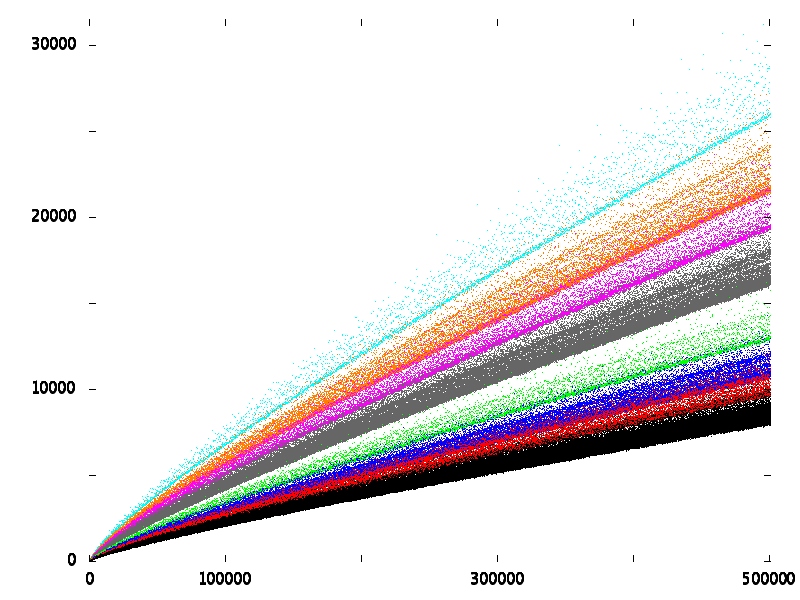} \\
\footnotesize Figura 6. Grafico di $\, g(n), \,$ per \vspace{-1.3mm} $\, n \in [\, 3 \, , \, 500.000\, ] \,$ \\
\hspace*{-10mm} Rispettivamente in nero, rosso, blu, \vspace{-1.5mm} verde, \\
\hspace*{-25.5mm} grigio, violetto, arancione, \vspace{-1.4mm} ciano \\
\hspace*{4.5mm} i punti con \vspace{18mm} $\, (2n,210) = 2, \, 14, \, 10, \, 70, \, 6, \, 42, \, 30, \, 210.$ \\
\normalsize
\end{center}
\textbf{\textit{iv}}) \ Per $\, m = 2310, \,$ vi sono 480 classi prime e 16 possibili divisori $\, d \,$ \\
per le classi pari e si ha:
\begin{center}
\begin{tabular}{c|cccccccccc}
$d$ & 2 & 6 & 10 & 14 & 22 & 30 & 42 & 66 & 70 & 110 \\
\hline
$\sigma_m(d)$ & 135 & 270 & \underline{180} & 162 & 150 & \underline{\underline{360}} & 324 & 300 & 216 & 200
\end{tabular}
\end{center}
\begin{center}
\hspace*{-41.3mm} \begin{tabular}{c|cccccc}
$d$ & 154 & 210 & 330 & 462 & 770 & 2310 \\
\hline
$\sigma_m(d)$ & \underline{180} & 432 & 400 & \underline{\underline{360}} & 240 & 480 
\end{tabular}
\end{center}
Si nota che $\, \sigma_{2310}(10) = \sigma_{2310}(154) = 180 \,$ e $\, \sigma_{2310}(30) = \sigma_{2310}(462) = 360; \,$ cos\`{\i} i divisori pari di $\, m = 2310 \,$ sono\, 16,\, ma le strisce che in questo caso si possono distinguere con i colori nel grafico della $\, g(n) \,$ sono di meno, solo 14 (fig.$\, 7$). Ovviamente per tutti i moduli $\, m = 2 \cdot 3 \cdots p_h \,$ con $\, h \geq 5 \,$ la $\, \sigma_m \,$ ristretta all'insieme dei divisori positivi pari di $\, m \,$ \underline{non} pu\`o essere iniettiva.
\begin{center}
\includegraphics[scale =.4]{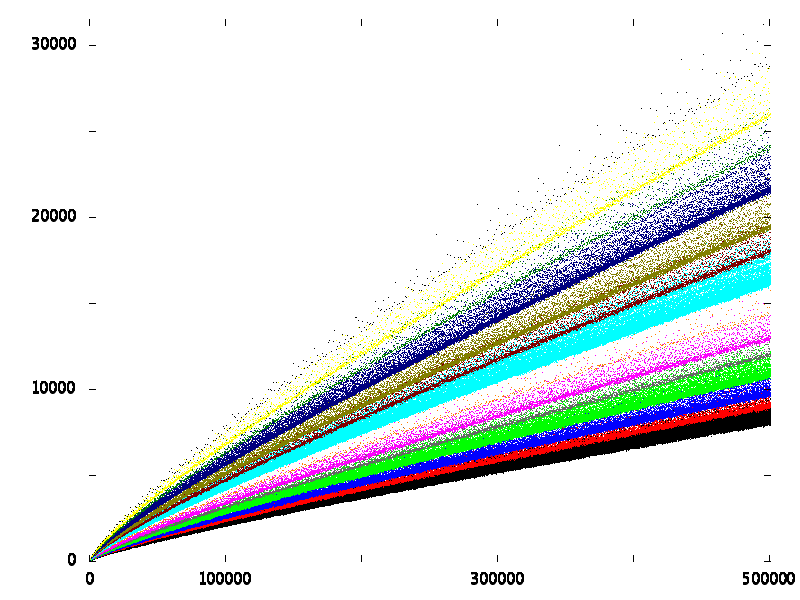} \\
\footnotesize Figura 7. Grafico di $\, g(n), \,$ per \vspace{-1.mm} $\, n \in [\, 3 \, , \, 500.000\, ] \,$ \\
\hspace*{-16mm} 14 colori, $\quad m = 2 \cdot 3 \cdot 5 \cdot 7 \cdot 11 = 2310. \,$ \\
\normalsize
\end{center}
\vspace{3mm}
\paragraph{3 Le altre.}
Un'analisi ampia e mirata, anche se di natura prevalentemente empirica, ha suscitato in chi scrive la convinzione che la validit\`a della congettura di Goldbach sarebbe eventualmente causata solo dal ``\textit{tipo di collocazione dei primi fra i naturali}'' e da nessun'altra loro propriet\`a intrinseca. A sostegno di questa convinzione vengono descritte alcune successioni di numeri dispari, sostanzialmente diverse dalla successione dei primi, ma che in un certo senso, presentano lo stesso ``\textit{tipo}\,'' di distribuzione fra i naturali e naturalmente la potenzialit\`a di soddisfare l'enunciato analogo alla congettura di Goldbach. \quad \`E necessario premettere qualche definizione. \\
Data una successione di numeri dispari $\, \mathfrak{r} = \{r_i \}_{i \in \mathbb{N}} \,$,\, monotona crescente, \\
si indicher\`a con $\, g_{\mathfrak{r}}(n) \,$ e si dir\`a ``\textit{funzione di Goldbach correlata alla successione $\mathfrak{r}$},'' la funzione aritmetica  che associa ad $\, n \,$ il numero delle coppie $\, (r_i,r_j) \,$ tali che $\, r_i+r_j=2n, \,$ \\
si dir\`a poi ``\textit{analogo della congettura di Goldbach per la successione $\mathfrak{r}$}\,'' \\
l'enunciato: ``\textit{vi \`e un} $\, m_\mathfrak{r} \in \mathbb{N} \,$ \textit{tale che} $\, g_{\mathfrak{r}}(n) > 0 \,$ \textit{per ogni} $\, n > m_{\mathfrak{r}} \,$''. \\
Si indicher\`a infine con $\, \pi_{\mathfrak{r}}(x) \,$ il numero degli elementi della successione $\, \mathfrak{r} \,$ che non superano $\, x.$ \\
\textbf{a}) \ Si consideri la successione di numeri dispari $\, \mathfrak{r} = \{r_i \}_{i \in \mathbb{N}} \,$ definita \\
come segue:
$$\left\{
\begin{array}{ll}
r_1 = p_2 & (p_2 = 3), \\
r_i = p_i+2 & \textrm{per} \ \; i > 1,
\end{array}
\right. $$
\`e facile vedere che $\, g(n) \leq g_{\mathfrak{r}}(n + 2) \leq g(n) + 2 \ \ \textrm{per ogni} \ \ n > 2; \,$ \\
ma gli elementi della $\, \mathfrak{r} \,$ sono ``\textit{quasi}\,'' tutti composti (per una ``buona met\`a'' sono addirittura multipli di 3 e i primi che vi figurano sono le seconde componenti delle coppie di primi gemelli). Per $\, N \,$ sufficientemente grande, \\
il grafico della $\, g_{\mathfrak{r}}(n) \,$ in $\, [3,N], \,$ sebbene diverso, sembra proprio la cometa di Goldbach. \\
Una successione analoga alla precedente, poco pi\`u complicata ma meno eccentrica, $\, \mathfrak{s} = \{s_i \}_{i \in \mathbb{N}} \,$,\, si ottiene ponendo:
$$\left\{
\begin{array}{l}
s_1 = p_2, \\
s_2 = p_3, \\
\hspace{30mm} \textrm{e \ per \ } i>2, \\
s_i = \left\{
\begin{array}{ll}
p_i+2, & \textrm{se} \ \; p_i = 6h-1, \\
p_i+4, & \textrm{se} \ \; p_i = 6k+1. \\
\end{array}
\right.
\end{array}
\right. $$
In questa successione compaiono ``\textit{in egual misura}\,'' dispari della forma $\, 6j \mp 1; \,$ vi sono inoltre infiniti numeri composti,\! \footnote{ \ \ Questa affermazione segue facilmente dal fatto che \ $\, \pi_{_{6,5}}(x) \sim \frac{x}{2 \lg x}, \,$ \ per $\, x \rightarrow \infty$ \\
\hspace*{7.5mm} e dal teorema di Brun: \ ``\textit{se} $\, \pi_{_2}(x) \,$ \textit{indica il numero delle coppie} $\, (p,p+2)$ \\
\hspace*{7.5mm} \textit{di primi gemelli con} $\, p+2 \leq x, \,$ \textit{vi \`e un intero} $\, n_{_0}, \,$ \textit{effettivamente computabile}, \\
\hspace*{7.5mm} \textit{tale che per} $\,x > n_{_0}, \,$ \textit{\`e} \ $\, \pi_{_{2}}(x) < \frac{100x}{(\lg x)^2} \,$'' \ [B].}\, nessuno dei quali divisibile per 3. \\
Anche in questo caso, il grafico della funzione correlata $\, g_{\mathfrak{s}}(n) \,$ in $\, [3,N], \,$ \\
con $\, N \,$ sufficientemente grande, figura come una cometa apparentemente indistinguibile da quella di \vspace{1.5mm} Goldbach. \\
\textbf{b}) \ La successione $\, \mathfrak{t} = \{t_h \}_{h \in \mathbb{N}} \,$,\, dove, posto $\, a_h = (h+1) \lg (h+1),$
$$t_h=\left\{
\begin{array}{ll}
\lfloor a_{h+1} \rfloor, & \textrm{se} \ \; \lfloor a_{h+1} \rfloor \ \; \textrm{\`e dispari,} \\
\lfloor a_{h+1} \rfloor + 1, \qquad & \textrm{altrimenti;}
\end{array}
\right. $$
trae la propria origine dalla \textit{formula asintotica dell'ennesimo primo}: \footnote{ \ \ Questo risultato \`e un enunciato equivalente al gi\`a citato ``\textit{teorema dei \vspace*{.5mm} numeri primi}.''}
$$\lim_{n\rightarrow\infty}\dfrac{p_n}{n\lg n} = 1.$$
Vi sono buoni motivi per credere che la successione $\, \mathfrak{t} \,$ soddisfi l'analogo della congettura di Goldbach e che $\, g_{\mathfrak{t}}(n) > 0 \,$ per ogni $\, n > 2.$
\begin{center}
\includegraphics[scale =.4]{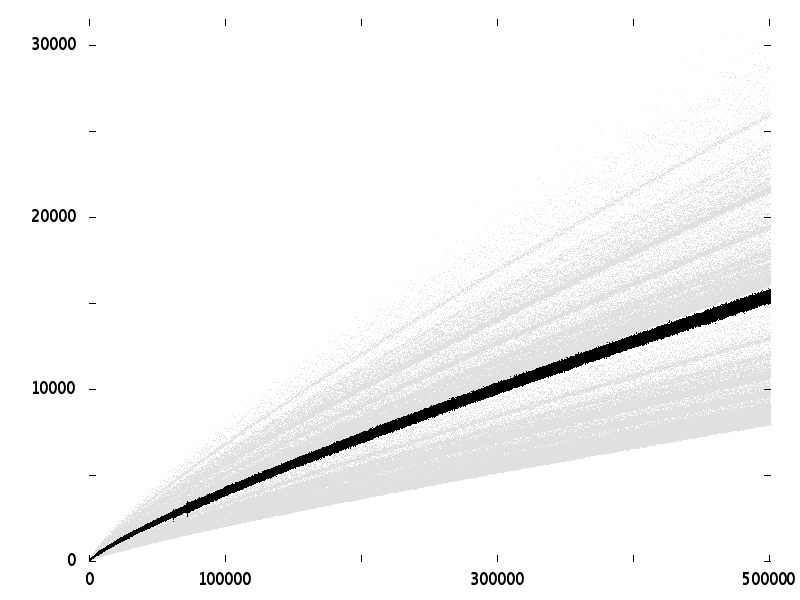} \\
\footnotesize Figura 8. Grafico di $\, g_{\mathfrak{t}}(n), \,$ per \vspace{-1mm} $\, n \in [\, 3 \, , \, 500.000\, ]$ \\
\hspace*{8mm} (in grigio chiaro il grafico di $\, g(n)).$ \\
\normalsize
\end{center}
Il grafico della $\, g_{\mathfrak{t}}(n) \,$ nel solito intervallo $\, [3,500.000] \,$ (fig.$\, 8$), bench\'e di notevole ``spessore,'' non sembra certo una cometa, ma si va a piazzare proprio nel giusto mezzo fra le due zone principali della cometa di Goldbach. \\
Si potrebbe eccepire che gli elementi della successione $\, \mathfrak{t} \,$ sono, fra i naturali, pi\`u ``\textit{frequenti}\,'' dei numeri primi, nel senso che $\, \pi_{\mathfrak{t}}(x) > \pi(x), \,$ per $\, x \geq 115.$ \footnote{ \ \ Tuttavia una ``\textit{frequenza}\,'' elevata non basta ad assicurare la validit\`a \vspace{-.6mm} dell'analogo \hspace*{7.5mm} della congettura di Goldbach; la successione dei naturali congrui a 3 modulo 4 \vspace{-.6mm} ne \`e \hspace*{7.5mm} un esempio molto semplice e significativo.}\, Effettivamente \`e noto che $\, n \lg n \leq p_n \, , \,$ per $\, n \geq 2; \,$ ma \`e anche \\
$\, n \lg (n \lg n) \geq p_n \, , \,$ per $\, n \geq 6 \;$  [R].\, \footnote{ \ \ Ovviamente $\, p_n \sim n \lg (n \lg n), \,$ per \vspace{-.7mm} $\, n \rightarrow \infty. \,$ \\
\hspace*{7.5mm} Sull'ennesimo primo vi sono risultati pi\`u stringenti; vanno certamente \vspace{-.7mm} citati: \\
\hspace*{7.5mm} [Cm] e \vspace{.7mm} [Dp].} \\
Cos\`{\i} \vspace{-4mm} posto
$$y_k =  (k+1) \lg \big((k+1) \lg(k+1) \big), \hspace{30mm} (2)$$
$$\textrm{e} \hspace{10mm} u_k=\left\{
\begin{array}{ll}
\lfloor y_{k+1} \rfloor, & \ \; \textrm{se} \ \; \lfloor y_{k+1} \rfloor \ \textrm{\`e dispari,} \\
\lfloor y_{k+1} \rfloor + 1, \qquad & \ \; \textrm{altrimenti,}
\end{array}
\right. $$
\hspace*{18mm} per ogni \vspace{1mm} $\, k \in {\mathbb{N}}, \,$ \\
si ha una successione di interi dispari $\, \mathfrak{u} = \{u_k \}_{k \in \mathbb{N}} \,$,\, monotona crescente, con $\, \pi_{\mathfrak{u}}(x) < \pi(x) \,$ per $\, x \geq 2, \,$ che, sebbene presenti una funzione correlata $\, g_{\mathfrak{u}}(n) \,$ con valori decisamente modesti rispetto alla $\, g_{\mathfrak{t}}(n) \,$ dell'esempio precedente, sembra comunque avere sufficienti possibilit\`a di verificare l'analogo della congettura di Goldbach (fig$\, 9$).
\begin{center}
\includegraphics[scale =.4]{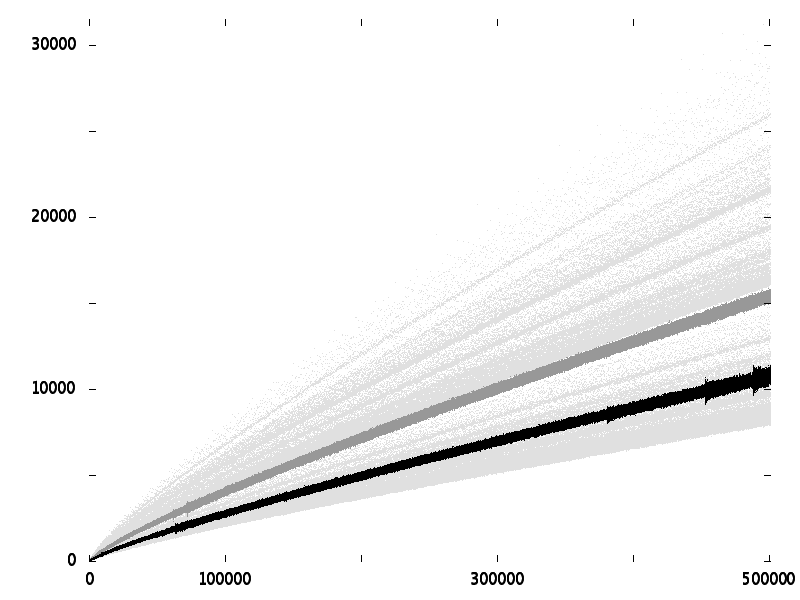} \\
\footnotesize Figura 9. Grafico di $\, g_{\mathfrak{u}}(n), \,$ per \vspace{-1mm} $\, n \in [\, 3 \, , \, 500.000\, ]$ \\
\hspace*{6mm} (rispett., in grigio e grigio \vspace{-1mm} chiaro \\
i grafici di $\, g_{\mathfrak{t}}(n) \,$ e $\, g(n)). \ $ \\
\normalsize
\end{center}
Vale la pena di notare che \`e $\, g_{\mathfrak{u}}(2l) = 0, \,$ per $\, 2 \leq l \leq 10 \,$ e pertanto \vspace{1.5mm} $\, m_{\mathfrak{u}} = 20. \,$  \\
\textbf{c}) \ Mediante la stessa successione (di numeri reali) $\, \{y_k \}_{k \in \mathbb{N}} \,$ definita sopra dalla\, (2),\, \`e possibile mostrare come vi sia un'ampia molteplicit\`a di successioni di numeri dispari, che pure si ritiene possano soddisfare l'analogo della congettura di Goldbach; tali sono le successioni di interi \underline{dis}p\underline{ari} $\, \mathfrak{v} = \{v_k \}_{k \in \mathbb{N}} \,$ con \\
\hspace*{30mm} $\displaystyle\left\{
\begin{array}{ll}
v_1 = 3 \hspace{9mm} & \textrm{e}\\
v_k \in [y_k, \, y_{k+1}) & \textrm{per} \ k \geq 2,
\end{array}
\right.$ \vspace{3mm} \\
tutte ovviamente monotone crescenti e con $\, \pi_{\mathfrak{v}}(x) < \pi(x) \,$ per $\, x \geq 2.$ \\
Sono state esaminate parecchie successioni di questa \textit{famiglia}, individuate sia effettuando sistematicamente la scelta di $\, v_k \,$ (ad esempio il massimo dispari in $\, [y_k, \, y_{k+1})), \,$  sia impiegando un generatore di numeri pseudocasuali; \\
in generale, i grafici delle funzioni $\, g_{\mathfrak{v}}(n) \,$ ad esse correlate, a meno di particolari minuti, rimangono assai simili al ``prototipo'' illustrato in figura \nolinebreak[4] 10.
\begin{center}
\includegraphics[scale =.4]{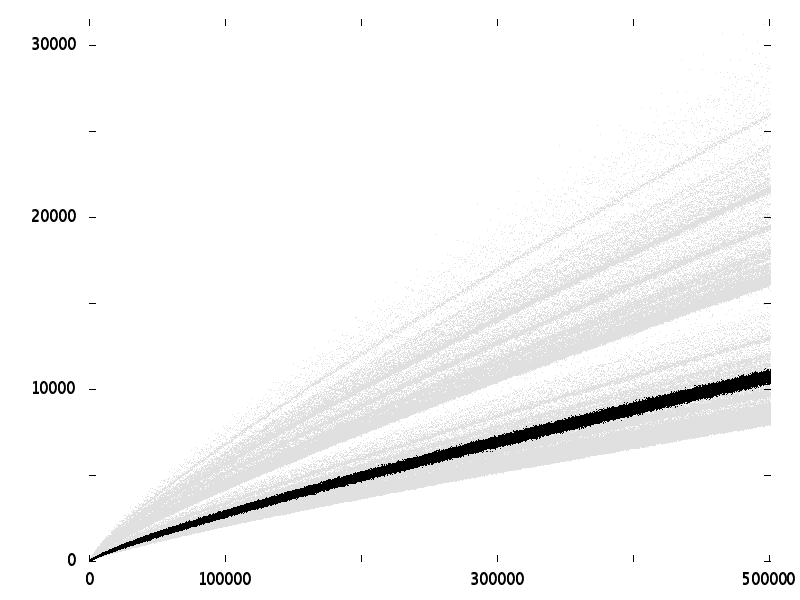} \\
\footnotesize Figura 10. Grafico di una $\, g_{\mathfrak{v}}(n), \,$ per \vspace{16mm} $\, n \in [\, 3 \, , \, 500.000\, ] \,$
\normalsize
\end{center}
\textbf{d}) \ Un'altra interessante famiglia di successioni crescenti di dispari \vspace{1.5mm} pu\`o essere individuata mediante la funzione $\, \dfrac{x}{\lg x - 1}; \,$ si vede facilmente che per ogni $\, j \in \mathbb{N}, \,$ l'equazione $\, \dfrac{x}{\lg x - 1} = 7 + j \,$ ammette un'unica soluzione \vspace{1.5mm} reale $\, x_j > 7; \,$ la successione $\, \{x_j \}_{j \in \mathbb{N}} \,$ di tali soluzioni risulta essere crescente e cos\`{\i} \`e possibile considerare le successioni $\, \mathfrak{w} = \{w_j \}_{j \in \mathbb{N}} \,$ i cui elementi sono interi \underline{dis}p\underline{ari} tali \vspace{3mm} che \\
\hspace*{30mm} $\displaystyle\left\{
\begin{array}{l}
w_1 = 3, \ \ w_2 = 5, \ \ w_3 = 7, \ \ w_4 = 11 \ \ \textrm{e} \\
w_{j+4} \in [x_j, \, x_{j+1}) \qquad \forall j \in \mathbb{N}.
\end{array}
\right.$ \vspace{3mm} \\
\`E noto che $\, \dfrac{x}{\lg x - 1} \leq \pi (x), \,$ per $\, x > 5393 \,$ (cfr. [Dp] pp. 37-41) e \vspace{1.5mm} risulta $\, \pi_{\mathfrak{w}}(x) \leq \pi (x), \,$ per $\, x \geq 2.$ \\
I grafici delle funzioni $\, g_{\mathfrak{w}}(n), \,$ correlate a diverse di queste successioni (determinate con scelta sistematica oppure pseudocasuale di $\, w_j), \,$ non presentano diversit\`a di rilievo dal campione illustrato in figura 11.
\begin{center}
\includegraphics[scale =.4]{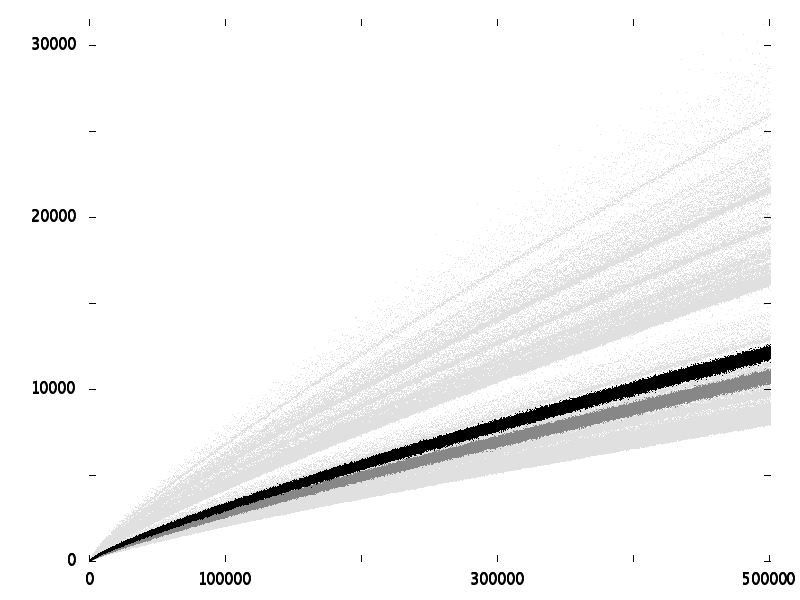} \\
\footnotesize Figura 11. Grafico di una $\, g_{\mathfrak{w}}(n), \,$ per $\, n \in [\, 3 \, , \, 500.000\, ] \,$ \\
\hspace*{.5mm} (in grigio scuro il grafico di una $\, g_{\mathfrak{v}}(n)).$ \\
\normalsize
\end{center}
Si pu\`o osservare che $\, \pi_{\mathfrak{w}}(x) \sim \dfrac{x}{\lg x} \sim \pi_{\mathfrak{v}}(x), \,$ quali che \vspace{1.5mm} siano $\, \mathfrak{w} \,$ e $\, \mathfrak{v} \,$ nelle rispettive famiglie e questo vale anche per le successioni $\, \mathfrak{r}, \ \,\mathfrak{s}, \ \, \mathfrak{t} \,$ e $\, \mathfrak{u}. \,$ L'idea di introdurre successioni $\, \mathfrak{c} \,$ di numeri naturali \vspace{1.5mm} con $\, \pi_{\mathfrak{c}}(x) \,$ asintoticamente equivalente a $\, \dfrac{x}{\lg x}, \,$ cio\`e con una distribuzione analoga a quella dei \vspace{1.5mm} numeri primi, sembra risalire almeno al 1937, ad opera di Harald Cram\'er, seppure con uno scopo diverso ([Ch], [Ga], [LZ]). \\
Qui si vuole solo porre in evidenza che a fianco alla congettura di Goldbach ve ne sono molte (infinite) altre analoghe. Francesco Lacava suggerisce che la validit\`a di una qualsiasi di queste congetture implicherebbe la validit\`a di tutte le altre. Sembra pi\`u verosimile che la risoluzione di qualche congettura analoga a quella di Goldbach, anche se non ne implica la risoluzione, potrebbe forse evidenziarne qualche nuova prospettiva.
\paragraph{4 Altre comete?}
Sembra che le funzioni $\, g_{\mathfrak{v}}(n) \,$ e $\, g_{\mathfrak{w}}(n) \,$ correlate alle successioni $\, \mathfrak{v} \,$ e $\, \mathfrak{w}, \,$ appartenenti alle due famiglie definite nel paragrafo precedente, abbiano tutte grafici somiglianti a ``scie'' puttosto che a ``comete'';  \`e possibile fra queste successioni sceglierne qualcuna che dia luogo ad una funzione che abbia il grafico simile ad una ``cometa''? La risposta \`e affermativa, basta preferire, nella scelta dell'elemento $\, w_j \,$ ($\, v_k \,$) di una successione $\, \mathfrak{w} \,$ ($\, \mathfrak{v} \,$) nell'intervallo $\, [x_j, \, x_{j+1}) \,$ ($\, [y_k, \, y_{k+1}) \,$), numeri (dispari possibilmente composti) che \underline{siano} p\underline{rimi} con $\, 3 \cdot 5 \cdots p_h \ \;  (h \geq 5). \,$ Perch\'e cos\`{\i} facendo gli elementi della successione $\, \mathfrak{w} \,$ ($\, \mathfrak{v} \,$) vengono distribuiti solo nelle classi prime col modulo $\, 2 \cdot 3 \cdots p_h. \,$
\begin{center}
\includegraphics[scale =.4]{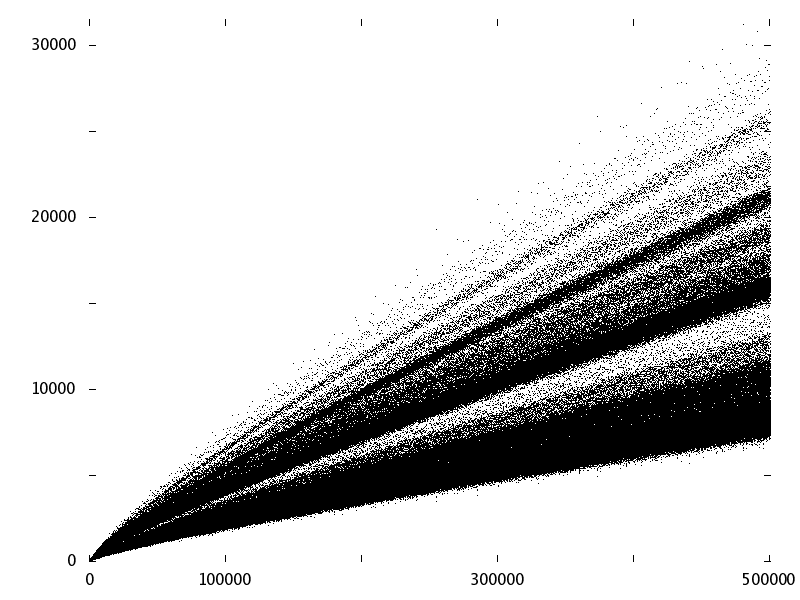} \\
\footnotesize Figura 12. Grafico di $\, g_{\mathfrak{w}'}(n), \,$ per $\, n \in [\, 6 \, , \, 500.000\, ] \,$
\normalsize
\end{center}
\begin{center}
\includegraphics[scale =.4]{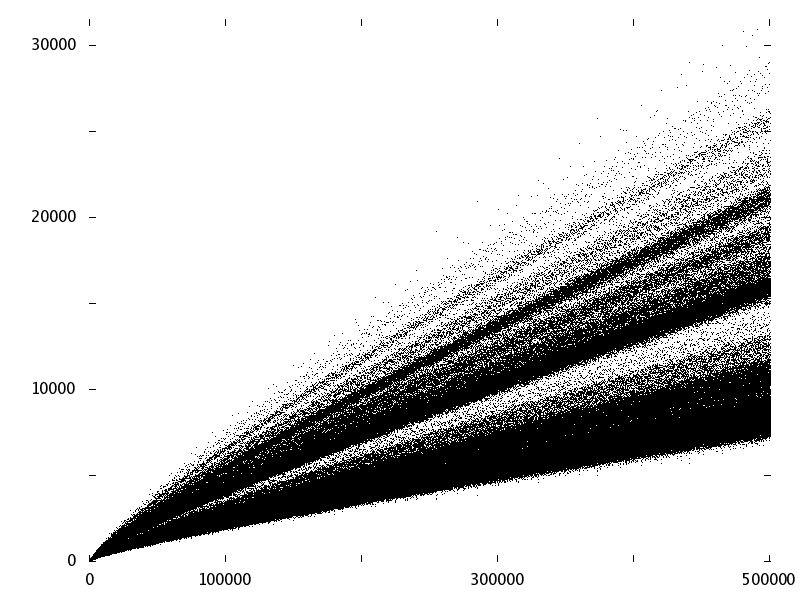} \\
\footnotesize Figura 13. Grafico di $\, g_{\mathfrak{w}''}(n), \,$ per $\, n \in [\, 6 \, , \, 500.000\, ] \,$
\normalsize
\end{center}
Le figure 12 e 13 mostrano i grafici di due funzioni $\, g_{\mathfrak{w}'}(n) \,$ e $\, g_{\mathfrak{w}''}(n) \,$ correlate con le successioni $\, \mathfrak{w}' \,$ e $\, \mathfrak{w}'' \,$ determinate ``preferendo'' dispari primi con $\, 3\cdot 5\cdot 7\cdot 11\cdot 13\cdot 17\cdot 19 = 4.849.845 \ \; (h =8); \,$ occorre notare che nella $\, \mathfrak{w}'' \,$ \`e stata accentuata al massimo la prevalenza di numeri \vspace{-1mm} composti.
\begin{center}
\small
Certo queste ``comete\,'' non sono belle quanto l'originale; \\
forse, qualcuno le direbbe ``\textit{artificiali},'' in realt\`a \vspace{16mm} esistono.
\end{center}
\large
\textbf{Bibliografia}\vspace{4mm} \\
\normalsize
\begin{tabular}{ll}
$[\textrm{B}]$ & \textsc{V. Brun}, Le crible d'Eratosth\`ene et le th\'eor\`eme de Goldbach. \\
& \textit{C.R. Acad. Sci. Paris} \textbf{168} (1919), 544-546. \\
$[\textrm{Bk}]$ & \textsc{K. Brown}, Evidence for Goldbach. \textit{Mathpages} (1994-2010). \\
& $http:\!//www.mathpages.com/home/kmath101.htm.$ \\
$[\textrm{Bn}]$ & \textsc{N. F. Benschop}, Additive structure of $\, \mathbb{Z}(\cdot)$ mod $m_k$ (square- \\
& free) and Goldbach's Conjecture. \textit{arXiv:math.GM}/0103091v5, \\
& 2009. $http:\!//arxiv.org/pdf/math/0103091$v5. \\
$[\textrm{C}]$ & \textsc{J. R. Chen}, On the representation of a large even integer \\
& as the sum of a prime and a product of at most two primes. \\
& \textit{Sci. Sinica} \textbf{16} (1973), 157-176. \\
$[\textrm{Ch}]$ & \textsc{H. Cram\'er}, On the order of magnitude of the difference \\
& between consecutive prime numbers. \textit{Acta Arith}. \textbf{2} (1937), 23-46. \\
$[\textrm{Cm}]$ & \textsc{M. Cipolla}, La determinazione assintotica dell'n$^{\textrm{imo}}$ numero \\
& primo. \textit{Matematiche Napoli} \textbf{3} (1902), 132-166. \\
$[\textrm{D}]$ & \textsc{L. E. Dickson}, \textit{History of the Theory of Numbers}. \\
& Vol. I, Dover Pub., (2005) 1919. \\
$[\textrm{Dm}]$ & \textsc{M. Deaconescu}, Adding units mod n. \\
& \textit{Elem. Math}. \textbf{55} (2000), 123-127. \\
$[\textrm{Dp}]$ & \textsc{P. Dusart}, Autour de la fonction qui compte le nombre \\
& de nombres premiers. \textit{Th\`ese de Doctorat - Universit\'e de Limoges} \\
& \textbf{17} 1998. $http:\!//www.unilim.fr/laco/theses/1998/T1998\underline{ }01.pdf.$ \\
$[\textrm{E}]$ & \textsc{T. Estermann}, On Goldbach's problem: proof that almost \\
& all even positive integers are sums of two primes, \\
& \textit{Proc. London Math. Soc.} (2) \textbf{44} (1938), 307-314. \\
\end{tabular} \\
\begin{tabular}{ll}
$[\textrm{Ga}]$ & \textsc{A. Granville}, Harald Cram\'er and the distribution \\
& of prime numbers. \textit{Scandanavian Actuarial J}. \textbf{1} (1995), 12-28. \\
& $http:\!//www.dms.umontreal.ca/~andrew/PDF/cramer.pdf$. \\
$[\textrm{Gr}]$ & \textsc{R. K. Guy}, \textit{Unsolved problems in number theory}. \\
& Springer-Verlag, 1994. \\
$[\textrm{HL}]$ & \textsc{G. H. Hardy - J. E. Littlewood}, \ Some problems \\
& of `partitio numerorum' III: on the expression of a number \\
& as a sum of primes. \textit{Acta Math.} \textbf{44} (1922) 1-70. \\
$[\textrm{LZ}]$ & \textsc{A. Languasco - A. Zaccagnini}, Intervalli fra numeri primi \\
& consecutivi. $http:\!//matematica$-$old.unibocconi.it/LangZac/$ \\
& \textit{home3.pdf}. \\
$[\textrm{MV}]$ & \textsc{H. L. Montgomery - R. C. Vaughan}, The exceptional set \\
& in Goldbach's problem. \textit{Acta Arith.} \textbf{27} (1975) 353-370. \\
$[\textrm{R}]$ & \textsc{J. B. Rosser}, Explicit Bounds for some functions of prime \\
& numbers. \textit{Amer. J. Math.} \textbf{63} (1941), 211-232. \\
$[\textrm{Rp}]$ & \textsc{P. Ribenboim}, \textit{The Little Book of Big Primes}. \\
& Springer-Verlag, 1991. \\
$[\textrm{S}]$ & \textsc{J. J. Sylvester}, On the partition of an even number into \\
& two primes. \textit{Proc. London Math. Soc. ser.}\;\raisebox{.5mm}{$_I$}\, \textbf{4} (1871), 4-6. \\
$[\textrm{Sh}]$ & \textsc{H. N. Shapiro}, \textit{Introduction to the Theory of Numbers}. \\
& Dover Pub., (2008) 1983. \\
$[\textrm{Sj}]$ & \textsc{J. W. Sander}, On the addition of units and nonunits mod m. \\
& \textit{Journal of Number Theory} \textbf{129} (2009), 2260-2266. \\
$[\textrm{V}]$ & \textsc{I. M. Vinogradov}, Representation of an odd number as a sum \\
& of three primes. \textit{Dokl. Akad. Nauk SSSR} \textbf{15} 6-7 (1937), 291-294. \\
$[\textrm{Z}]$ & \textsc{A. Zaccagnini}, Variazioni Goldbach: problemi con numeri \\
& primi. \textit{L'Educazione Matematica, Anno XXI Serie VI} \textbf{2} (2000),\\
&  47-57. $http:\!//www.math.unipr.it/\!\!\sim\!zaccagni/psfiles/papers/$ \\
& \textit{Goldbach\underline{ }I.pdf}.
\end{tabular}
\end{document}